\documentclass{amsart}

\usepackage{amsmath,amssymb,amsthm}
\usepackage{graphicx,tikz}

\newtheorem*{thm}{Theorem}

\newtheorem{corollary}{Corollary}
\newtheorem*{conjecture}{Conjecture}
\newtheorem{lemma}{Lemma}

\theoremstyle{definition}

\theoremstyle{remark}

\newcommand{\norm}[1]{\left\lVert#1\right\rVert}

\begin{document}

\title[]{Spectral Clustering revisited: Information hidden in the Fiedler vector}

\author[]{Adela DePavia}
\address{Yale University, New Haven, CT 06511, USA}
\email{adela.depavia@yale.edu}

\author[]{Stefan Steinerberger}
\address{Department of Mathematics, Yale University, New Haven, CT 06511, USA}
\email{stefan.steinerberger@yale.edu}
\thanks{AD was supported by the Yale University ESI PREP Post-Baccalaureate Research Education Program. SS was partially supported by the NSF (DMS-1763179) and the Alfred P. Sloan Foundation.}

\begin{abstract} We are interested in the clustering problem on graphs: it is known that if there are two underlying clusters, then the signs of the eigenvector corresponding to the second largest eigenvalue of the adjacency matrix  can reliably reconstruct the two clusters. We argue that the vertices for which the eigenvector has the largest and the smallest entries, respectively, are unusually strongly connected to their own cluster and more reliably classified than the rest. This can be
regarded as a discrete version of the Hot Spots conjecture and should be useful in applications. We give a rigorous proof for the stochastic block model and several examples.
\end{abstract}

\maketitle

\section{Introduction}
\subsection{Introduction.}

The purpose of this paper is to discuss a general refinement of the spectral clustering approach that seems very widely applicable. For simplicity of exposition, we will restrict
ourselves to the simplest possible case: suppose we are given a connected, unweighted, undirected graph $G$ that decouples into two equally sized clusters that
have many edges within themselves and very few edges across (see Fig. 1 for a sketch of what this could look like).

\begin{center}
\begin{figure}[h!]
\includegraphics[width=0.6\textwidth]{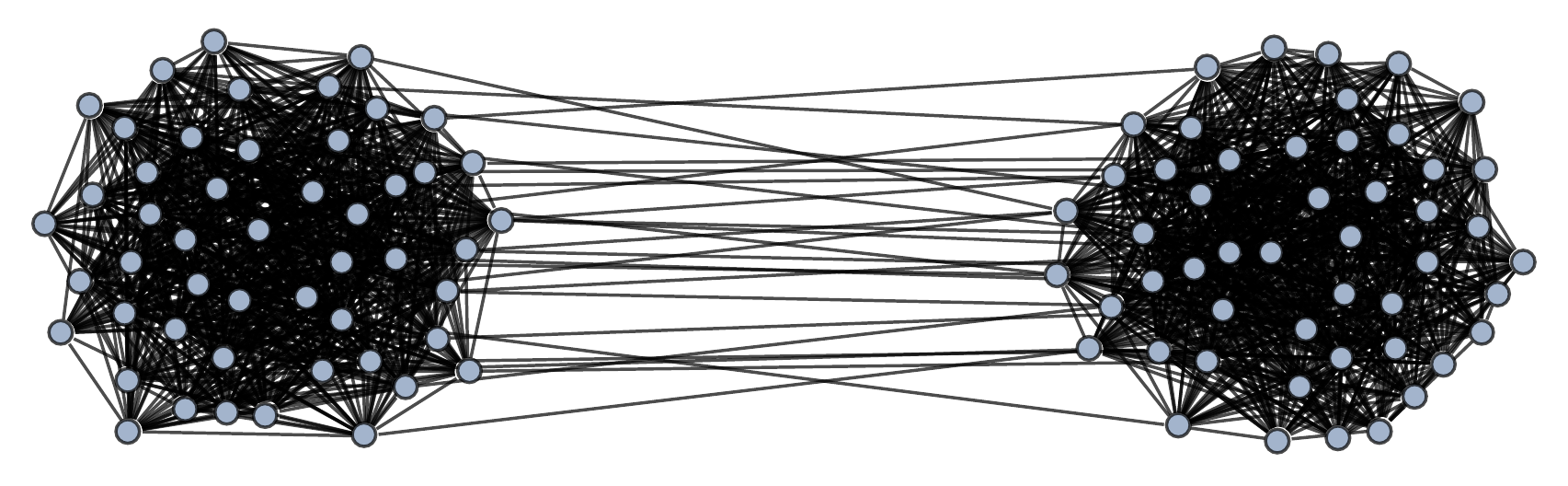}
\caption{A graph with two communities that have strong inter-connectedness but very few edges between them.}
\end{figure}
\end{center}

The spectral clustering approach is quite simple: let $A \in \mathbb{R}^{n \times n}$ be the adjacency matrix associated to the graph $G$. The matrix $A$ is symmetric and has real
eigenvalues and eigenvectors. It is known that if the underlying graph nicely decomposes into two roughly equally sized clusters with few connections between them, then the
second-largest eigenvector of $A$, we shall denote it by $\mathbf{v_2}$, is essentially constant on each cluster and, in particular, the signs of $\mathbf{v_2}$ allow us to reconstruct to which cluster
any specific vertex belongs. Nowadays, this is considered a quite classical construction and it is very well understood. We refer to \cite{blum, lux, ng, spiel, vershynin}.

\subsection{The Main Idea.}\label{ssec::main_idea} Our main idea is to not only look at the sign of the second eigenvector but also at the \textit{size} of the entry.

\begin{quote}
\textit{An Informal Spectral Clustering Hot Spots Conjecture.} Suppose we are using the sign of the second eigenvector $\mathbf{v_2}$ to partition a graph into two parts, the $(+)$ part and the $(-)$ part. If $i,j$ are two vertices in the graph and $\mathbf{v_2}(i) > \mathbf{v_2}(j) > 0$, then $i$ is `more likely' to truly be correctly identified as being in the $(+)$ cluster than $j$ (and likewise for the negative entries). 
\end{quote}

There are many ways of making this precise. Here is one natural (informal) conjecture: suppose $G$ is made up of two clusters of roughly equal size. Then the classification error on the extremal set
$$ E_{\delta}= \left\{i \in V:  \mathbf{v_2}(i)>0 ~\mbox{and} ~ \#\left\{j \in V: \mathbf{v_2}(j) > \mathbf{v_2}(i) > 0\right\} \leq \delta |V|  \right\}$$
is much smaller than the overall classification error. We prove this for the stochastic block model and show in \S 3 that this can be empirically observed. We point out that, while the idea is exceedingly natural, we are not aware of many theoretical results in this direction. Indeed, an analogous question in the continuous setting has been open since 1974
and is suspected to be quite difficult (see \S \ref{ssec::hot_spots} for a discussion of the Hot Spots conjecture). \\

 This raises an interesting question: if it is indeed the case that vertices corresponding to extremal values
of the eigenvector are more likely to be correctly identified by spectral clustering, is it then possible to propagate this high-quality classification from the `extremal vertices' to the rest?

\begin{quote}
\textit{An Algorithmic Challenge.} Is it possible to make use of the fact that the error rate for `extremal' entries is much smaller to improve on the overall classification error?
\end{quote}

\begin{center}
\begin{figure}[h!]
\begin{tikzpicture}
\node at (0,0) {\includegraphics[width=0.9\textwidth]{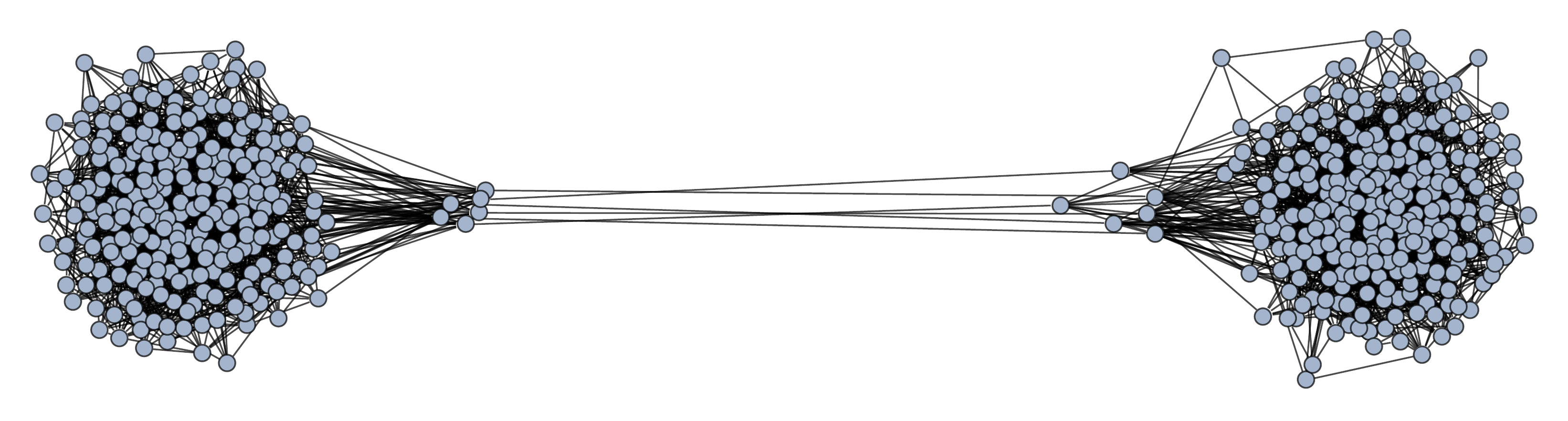}};
\node at (-2.2,1.59) {the eigenvector here};
\draw [thick, ->] (-1, 1.3) -- (-2,0.4);
\node at (1.93, 1.59 ) {is probably closer to 0 than here};
\draw [thick, ->] (5, 1.5) to[out=40, in=40] (5.6, 0.5);
\end{tikzpicture}
\caption{A stochastic block graph: $n=500$, $p=0.05$ and $q=0.0001$. Two clusters that are barely connected; the cluster on the left has a somewhat more `ambiguous' region; this ambiguity should be reflected in the \textit{size} of the entries of the eigenvector $\mathbf{v_2}$.}
\end{figure}
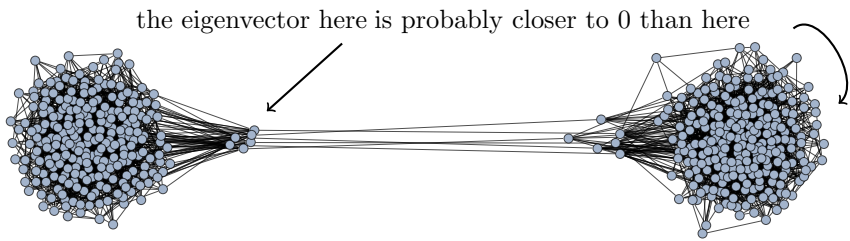
\end{center}

Our main contribution will be to illustrate the principle in the case of a stochastic block model (see Fig. 1 or Fig. 2). The stochastic block model is a model of random graphs. Let $n \in \mathbb{N}$ be an even integer and let $0 < q < p < 1$. We denote by $G(n,p,q)$ a type of random graph that generalizes the Erd\H{o}s-Renyi random graph: we assemble $n/2$ vertices in one group and the remaining $n/2$ vertices in another. We now consider all pairs of different vertices and add an edge connecting them with likelihood $p$ if they are in the same group and with likelihood $q$ if they are in different groups.
We will work in the case where $p > q$ are fixed and $n \rightarrow \infty$. This case is understood and it is known that the second eigenvector will identify all vertices correctly with high probability. Much stronger results (where $p$ and $q$ get closer to each other as $n$ increases) are known \cite{abbe2, abbe1, band, mc, rohe, vershynin}, see especially the recent survey \cite{abbe}.
One could expect the asymptotic behavior to be of the type
$$ \mathbf{v_2}(i) = \begin{cases} +1/\sqrt{n} \qquad \mbox{if}~i~\mbox{in the first cluster} \\ -1/\sqrt{n} \qquad \mbox{if}~i~\mbox{in the second cluster} \end{cases} + \mathcal{O}\left(\frac{\log{n}}{n}\right).$$
Speaking of the vertex for which the eigenvector $\mathbf{v_2}$ has the `largest' entry seems like a misnomer since the largest (positive) entry and the smallest (positive) entry of $\mathbf{v_2}$ are basically identical in size (their ratio tends to 1 as $n \rightarrow \infty$).
We will show that even in this rather degenerate case (wherein the underlying graph cannot be said to approximate any nice smooth manifold; indeed, the diameter is 2 with high probability), our main idea is still valid: the error term contains \textit{a lot} of information and the vertex with the largest entry is indeed
different from other vertices.

\subsection{The Hot Spots conjecture}\label{ssec::hot_spots}

Let us consider the continuous setting. Let $(M,g)$ be a smooth, compact manifold. A natural question is whether it is possible for the Rayleigh-Ritz quotient
\begin{align*} \label{eig}
 \lambda_2 = \inf_{f \neq 0} \frac{\int_{M} |\nabla f|^2 dg}{ \int_{M} f^2 dg}
 \end{align*}
 
 to be small. The volume $\mbox{vol}(M)$ plays a role. However, there is a second type of constraint (the discovery of which is attributed to Calabi \cite{cheeger}): for a manifold very similar to the one shown in Figure 3 (a `dumbbell' domain), it is possible to make the function $f$ essentially constant
on each side and make it smoothly varying from one component to the other thereby concentrating all the change in a tiny area. 
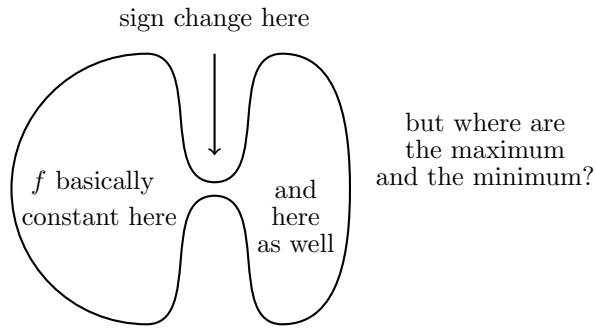
\begin{figure}[h!]
\begin{center}
\begin{tikzpicture}[scale=0.9]
\draw [thick] (0,0) to [out=270,in=180] (2,-2) to [out=0,in=180] (3,-0.1) to [out=0,in=180] (4,-2) to [out=0, in=270] (5,0) 
to [out=90, in = 0] (4,2) to [out=180, in = 0] (4,2) to [out=180, in = 0] (3,0.1) to [out=180, in = 0] (2,2) to [out=180, in = 90] (0,0);
\node at (1.2,0.1) {$f$ basically} ;
\node at (1.3,-0.4) {constant here} ;
\node at (4.2, 0) {and};
\node at (4.2, -0.4) {here};
\node at (4.2, -0.8) {as well};
\node at (3, 2.5) {sign change here};
\draw [thick, ->] (3,2) -- (3,0.5);
\node at (7,1) {but where are};
\node at (7,0.6) {the maximum};
\node at (7,0.2) {and the minimum?};
\end{tikzpicture}
\end{center}
\caption{A two-dimensional domain with Neumann boundary condition and the second Laplacian eigenfunction.}
\end{figure}

The gradient $\nabla f$ would not be small but it only plays a role in the
`tube' connecting the two balls and we can make this tube have arbitrarily small area by making it thinner.
Cheeger, in establishing the celebrated Cheeger inequality \cite{cheeger}, showed that this is essentially the only obstruction: unless a structure
of this type arises, $\lambda_2$ cannot be arbitrarily small. Conversely, in many cases the function attaining $\lambda_2$ has exactly this structure. Moreover,
none of these ideas are restricted to manifolds and they survive the transition to the discrete graph setting. Cheeger inequalities have become a standard ingredient
in spectral graph theory, see e.g. \cite{chung, levin}.
Let us now assume the manifold is actually a nice domain $\Omega \subset \mathbb{R}^2$ with boundary and let $f:\Omega \rightarrow \mathbb{R}$ be the function minimizing 
the Rayleigh-Ritz quotient among all functions $f:\Omega \rightarrow \mathbb{R}$ that have mean value 0 on $\Omega$ (see Fig. 2). This is known as the second eigenfunction of the (Neumann-)Laplacian and describes the generic long-time behavior of the heat equation in an insulated room the shape of $\Omega$.  A natural question is the following: where does $f$ assume its maximum and minimum? Going back to physical intuition, it seems reasonable to assume that maximum and minimum should be attained at the boundary. 

\begin{quote}\textbf{Hot Spots Conjecture} (strong form, $d=2$)\textbf{.} Let $\Omega \subset \mathbb{R}^2$ be simply connected. Then the maximum and the minimum of the second eigenfunction are located at the boundary.
\end{quote}

The Hot Spots conjecture dates back to a 1974 lecture of Rauch that he gave at a Tulane University PDE conference \cite{b3}. It was mentioned in a 1985 book of Kawohl \cite{kawohl} who also suggested that it may be false in general but could be true for convex domain.  Ba\~nuelos \& Burdzy \cite{b3} proved it for obtuse triangles and Burdzy \& Werner \cite{b4} obtained a counterexample for domains that are not simply connected (see also \cite{b25}). Judge \& Mondal \cite{judge} recently established the Hot Spots conjecture for all triangles; the second author proved that in a convex domain the maxima and minima are at a distance close to the diameter from each other \cite{steini}. The Hot Spots conjecture is now widely assumed to be true for convex domains (and possibly even for simply connected domains). It cannot be easily translated to the discrete setting: after all, what is the boundary of a graph? 
\begin{quote}
\textit{Another Informal Spectral Clustering Hot Spots Conjecture.} Suppose we are using the sign of the second eigenvector $\mathbf{v_2}$ to partition a graph into two parts, the $(+)$ part and the $(-)$ part. Large entries of $\mathbf{v_2}$ correspond to vertices having a disproportionately large number of neighbors in the same cluster (and similarly for the smallest entries).
\end{quote}
However, various reformulations are meaningful: in the case of the stochastic block model, a particularly canonical formulation is the main result of this paper. It would be interesting to see what kind of results are possible for other types of random graphs or possibly even in the deterministic setting.

\subsection{Related results.}
We are not aware of this conjecture and this challenge being stated anywhere explicitly. Based on work of Rachh and the second author \cite{manas} in the continuous case, Cheng, Rachh and the second author \cite{xiu} proved that the maximum and the minimum on graphs are attained `far away' from the vertices adjacent to a sign change. A similar idea was used by Cheng, Mishne and the second author for averages over eigenfunctions \cite{gal}. However, generally not much seems to be known (and the difficulty of the Hot Spots conjecture in the continuous setting may partially explain why). It is known, under very general conditions, that the sign of the second eigenvector splits the graph into two connected graphs (one where the eigenvector is positive and one where it is negative) -- this is the seminal work of Fiedler \cite{fiedler1, fiedler2, fiedler3}. However, apart from those results, very little is known about how the Fiedler vector is connected to topological properties of the underlying graph \cite{gern}.
The Hot Spots conjecture has been stated explicitly on graphs.  Chung, Seo, Adluru \& Vorperian \cite{moo} asked whether maximum and minimum are attained at vertices at maximal distance of one another; this is false in general but sometimes true \cite{roy}.

\section{Results}

\subsection{Stochastic Block Model.} For even integers $n \in \mathbb{N}$ and given parameters $0 < q < p < 1$, we denote by $G(n,p,q)$ a type of random graph that is constructed as follows: we put the first $n/2$ vertices in one group and the remaining $n/2$ vertices in another. We now consider all pairs of different vertices and add an edge connecting them with likelihood $p$ if they are in the same group and with likelihood $q$ if they are in different groups. 
It is known (we will recall the argument below) that the first eigenvalue of the adjacency matrix $A$ satisfies 
$$\lambda_1 = \frac{(p+q)n}{2} + \mathcal{O}(\sqrt{n}) \qquad \mbox{with high probability.}$$
Moreover, the associated eigenvector $\mathbf{v_1}$ is close to the constant vector.
The second eigenvalue of the adjacency matrix is known to be
$$\lambda_2 = \frac{(p-q)n}{2} + \mathcal{O}(\sqrt{n}) \qquad \mbox{with high probability.}$$
Moreover, the second eigenvector $\mathbf{v_2}$ is close to a vector that is constant on each cluster (and has mean value close to 0). In particular, the sign of the entries of the second eigenvector $\mathbf{v_2}$ can be used to identify which vertex belongs to which group. These things are by now fairly classical; we refer to the very clear exposition in the textbook of Vershynin \cite[\S 4.5.]{vershynin} and references therein.

\subsection{The Main Result.} We are now ready to state the main result. We consider the classical case of $p>q$ (more in-group connections than out-group connections) and study the behavior of the eigenvector $\mathbf{v_2}$ as $n$ becomes large. We assume that $\mathbf{g} \in \left\{-1,1\right\}^n$ is the vector indicating group membership of each vertex. A classical approach to spectral clustering is based on some elementary facts about random matrices which we survey below in \S \ref{ssec::theorem_prelims} and which imply
$$ \mathbf{v_2} = \frac{\mathbf{g}}{\sqrt{n}} + \mathbf{e_2} \qquad \mbox{where}  \quad \| \mathbf{e_2} \| \lesssim  \frac{1}{\sqrt{n}} \frac{\| \textbf{g} \|}{\sqrt{n}} \quad \mbox{with high probability.}$$
The purpose of our paper is to show that the entries of $\mathbf{e_2}$ are actually highly meaningful and encode a lot of the underlying information as $n \rightarrow \infty$. In particular, it will allow us to deduce a statement implying that the vertex on which the vector $\mathbf{v_2}$ has its largest absolute-value entry is actually somewhat distinguished in the sense that it has disproportionately many neighbors that lie in the same cluster. This relationship between the magnitude of $\mathbf{v_2}$ and a node's attachment to its ``correct'' cluster is not only true for the largest magnitude entries: we show that the size of $\mathbf{v_2}$ encodes meaningful information about all nodes.

\begin{thm}[Hidden Information] Let $0 < q < p < 1$ and $\varepsilon_0 > 0$. Then, as $n \rightarrow \infty$
 \begin{align*}
  \mathbf{e_2} &=  \frac{2}{(p-q)n^{3/2}}  \left(A\mathbf{g} -  \lambda_2 \mathbf{g} \right)  + \mathbf{error},
  \end{align*}
 where $\| \mathbf{error} \| \lesssim n^{-1+\varepsilon_0}$ with high probability.
\end{thm}
It is known that $\|\mathbf{e_2}\| \lesssim n^{-1/2}$ with high probability (see \S \ref{ssec::theorem_prelims}). This shows that the Theorem does indeed
capture the entire expansion of $\mathbf{e_2}$ up to a lower order term at size $\lesssim n^{-1 + \varepsilon_0}$. 
The result is with high probability in the usual sense (the likelihood of it failing decays exponentially
in $n$; the proof would allow for this decay rate to be made quantitative but this is perhaps not quite as interesting). We now turn to the interpretation of the result. 
We like to think of the approximation as being comprised of two parts: by rewriting the algebraic expression and ignoring $\mathbf{error}$, we arrive at
$$ \mathbf{e_2} \sim  \underbrace{ \frac{2}{(p-q)n^{3/2}}   \left(A\mathbf{g} - \frac{(p-q)n}{2}\mathbf{g}\right)}_{\tiny \mbox{local property of vertices}} - \underbrace{\frac{2 \mathbf{g}}{(p-q)n^{3/2}}   \left( \lambda_2 - \frac{p-q}{2}n\right)}_{\tiny \mbox{global shift}}.$$
Both terms are vectors at the same scale $\sim n^{-1/2}$. We will now explain both.

\subsubsection{Local properties.}\label{ssec::local_shift}
Let $i$ correspond to a vertex in the first cluster (which we define without loss of generality to be the cluster where the indicator vector $\textbf{g}$ has positive entries). An interesting local quantity is how many of the neighbors of $i$
are in the same cluster and how many neighbors of $i$ are in the others cluster. This motivates the quantity
$$ d_i = [\#\mbox{neighbors of}~i~\mbox{in the same cluster}] -  [\#\mbox{neighbors of}~i~\mbox{in the other cluster}].$$
By the central limit theorem, we expect
$$ d_i = \frac{p-q}{2}n  \pm \mathcal{O}(\sqrt{n}).$$
What we are interested in is the error term: it is going to be of size $\sim \sqrt{n}$ but how big is it?  Is the sign positive or negative? 
The answer is given by the $i-$th entry of the vector $\Delta = (\Delta_1, \dots, \Delta_n)$, where $\Delta$ is given by
$$ \Delta = A\mathbf{g} - \frac{(p-q)n}{2}\mathbf{g} \qquad \mbox{which arises as the first term in our expansion}.$$
The same is true, up to a flip of the sign, for the vertices in the second cluster. In summary: the first term in our expansion of $\mathbf{e_2}$, encoding local property of the vertices, contains the deviation from the expected
number of `in-cluster neighbors' minus `out-cluster neighbors' (and with a flipped sign for the vertices where the indicator vector is negative). 
We expect each entry of $\Delta$ to be
well approximated by a Gaussian centered at 0 with standard deviation $\sim_{p,q} \sqrt{n}$. This shows that, with high probability, the first term in the expansion is of size $\sim n^{-1/2}$
as predicted.
This means that vertices with a disproportionately large number of neighbors within their own cluster have a \textit{slightly} larger (in absolute value) entry in $\mathbf{v_2}$.

\subsubsection{Global shift.}\label{ssec::global_shift} The global shift is quite easy to understand. Observe that the global shift term contains the deviation from the expected value of $\lambda_2$. It is known that
$$ \lambda_2 = \frac{p-q}{2} \pm \mathcal{O}(\sqrt{n}) \qquad \mbox{with high probability.}$$
This shows that we also expect the global shift to be of size $\sim n^{-1/2}$. Moreover, we recall that
$$ \mathbf{v_2} = \frac{\mathbf{g}}{\sqrt{n}} + \mathbf{e_2} $$
which means that the global shift, being a multiple of $\mathbf{g}$, can be absorbed in the first term.  Put differently, the global shift is actually constant on the first cluster and constant on the second cluster and therefore merely
shifts values but does not have any impact on which vertex has the largest entry or even the relative ordering among the entries. This can also be seen from rewriting $\mathbf{v_2}$ by moving
the global shift $\varepsilon_{\text{global}}$ into the leading term 
\begin{equation}\label{eq::v_2_expansion}
    \mathbf{v_2} = (1+\varepsilon_{\text{global}})\frac{\mathbf{g}}{\sqrt{n}}  + \frac{2}{(p-q)n^{3/2}}   \left(A\mathbf{g} - \frac{(p-q)n}{2}\mathbf{g}\right)  + \mathbf{error}.
\end{equation}

\subsubsection{Summary.} We have shown that $\mathbf{v_2}$ is essentially, up to an error at a smaller scale, given by $\mathbf{g}$ and a variation sitting on top that describes the number of in-cluster neighbors minus the number of out-cluster neighbors. We observe that this is, at the same time, the dominant form of randomness governing $\mathbf{v_2}$. We observe that, for any individual node,
$$ \mathbb{E} \#~\mbox{in-cluster-neighbors} = \frac{pn}{2}, \qquad  \mathbb{V} \#~\mbox{in-cluster-neighbors} =\frac{n p (1-p)}{2}.$$
$$ \mathbb{E} \#~\mbox{out-cluster-neighbors} = \frac{qn}{2}, \qquad  \mathbb{V} \#~\mbox{out-cluster-neighbors} =\frac{n q (1-q)}{2}$$
and thus each entry of this matrix satisfies (in a component-by-component sense)
$$  \mathbb{E} \left(A\mathbf{g} - \frac{(p-q)n}{2}\mathbf{g}\right) = \mathbf{0}$$
and
$$ \mathbb{V}  \left(A\mathbf{g} - \frac{(p-q)n}{2}\mathbf{g}\right) = n\left( \frac{p(1-p)}{2} + \frac{q(1-q)}{2} \right) \mathbf{1}.$$
Moreover, by the central limit theorem, the asymptotic behavior of these random variables starts to behave (rather quickly) like a Gaussian random variable at the same scales. Some obvious consequences are derived in the next section.

\subsection{Some Implications.}
We summarize our discussion until now. For the stochastic block model, the second eigenvector $\mathbf{v_2}$ is known to be asymptotically accurate in terms of cluster identification: the sign of each entry allows to
recover the cluster identity of the vertex. Moreover, under the ansatz
$$ \mathbf{v_2} = \frac{\mathbf{g}}{\sqrt{n}} + \mathbf{e_2},$$
then $\|\mathbf{e_2}\| \sim n^{-1/2}$. This means that we expect a typical entry of $\mathbf{e_2}$ to be of size $\sim n^{-1}$. The eigenvector $\mathbf{e_2}$ decouples into two terms, one
being a multiple of $\mathbf{g}$ that can be absorbed in the leading term, the other measuring whether a vertex has disproportionately more in-group neighbors than out-group neighbors. 

\begin{corollary}\label{corr::extreme_vertices_extremize_attachment} For every $\eta > 0$, there exists $\varepsilon >0$ such that at least $99\%$ the vertices corresponding to the the $\varepsilon\cdot n$ largest entries of $\mathbf{v_2}$ satisfy
$$ \# \emph{in-cluster neighbors} - \# \emph{out-cluster neighbors} \geq \eta \sqrt{n} \qquad \mbox{w.h.p.}$$
Likewise, at least $99\%$ of the verticies corresponding to the $\varepsilon\cdot n$ smallest entries of $\mathbf{v_2}$ will satisfy
$$ \# \emph{in-cluster neighbors} - \# \emph{out-cluster neighbors} \geq \eta \sqrt{n} \qquad \mbox{w.h.p.}$$
\end{corollary}

This argument can be made quantitative: $\varepsilon$ decreases extremely rapidly as $\eta$ increases. One could also slightly refine the $99\%$ and replace it by a quantity tending to 1 as $n \rightarrow \infty$ using the same argument. Given the precision of the asymptotic expansion, many other corollaries are conceivable.\\

 Corollary \ref{corr::extreme_vertices_extremize_attachment} also implies that the nodes with the most extreme-magnitude entries of $\mathbf{v_2}$ are more likely to be correctly classified by the sign of $\mathbf{v_2}$.  The corollary states that w.h.p. the nodes with the $\varepsilon\cdot n$ largest-magnitude entries of $\mathbf{v_2}$ also have a larger-than-expected 
\begin{align*} 
(\# \text{in-cluster neighbors} &- \# \text{out-cluster neighbors})\\
&- \mathbb{E}[\# \text{in-cluster neighbors} - \# \text{out-cluster neighbors}]
\end{align*} 
This is exactly equal to the local-shift term discussed in \S \ref{ssec::local_shift} and, by the main theorem, this difference is the dominant deviation to the value of $\mathbf{v_2}$ at each node. In particular, if $\# \text{in-cluster neighbors} > \# \text{out-cluster neighbors}$ for a given node (which occurs w.h.p. for extremal nodes in the regime where $p > q$), the sign of this additive deviation agrees with the sign of $\mathbf{g}$ at that node. Thus, Corollary \ref{corr::extreme_vertices_extremize_attachment} implies that w.h.p. $\mathbf{v_2}$-extremal vertices have the sign of their $\mathbf{v_2}$ entry agree with the sign of their indicator vector $\mathbf{g}$ entry, i.e. extremal vertices are more likely to be correctly classified, and thus error rates on extremal vertices are lower than the global error rate with high probability.\\

We conclude with a conjecture for the extremal behavior.
\begin{conjecture} The largest entry of $\mathbf{v_2}$ in the stochastic block model satisfies
$$ \| \mathbf{v_2} -   \frac{\mathbf{g}}{\sqrt{n}} \|_{\infty} \sim \frac{\log{n}}{n}.$$
\end{conjecture}
This is perhaps not surprising: given enough vertices, there is always one that has a disproportionate number of connections within their own cluster and this is then reflected in $\mathbf{e_2}$.
More precisely, we expect the difference between in-group neighbors and out-group neighbors to be dominated by a Gaussian (as the limiting object arising from the sum of many Bernoulli random variables).
It is known that for $n$ independent Gaussians $X_1, \dots, X_n$ where $X_i \sim \mathcal{N}(0,1)$, we have
$$ \mathbb{E} \max_{1 \leq i \leq n}{ X_i} \sim \log{n}$$
and this motivates the appearance of the logarithm in the conjecture.
The conjecture follows almost, but not quite, from our argument: we cannot, at the moment, exclude that
the \textbf{error} vector, which we know to be of size $n^{-1 + \varepsilon_0}$ concentrates entirely on the vertex where $\mathbf{v_2}$ assumes its maximum.

\section{Numerical Examples}
This section shows some numerical examples illustrating the main points of this paper. We start by discussing the ideal cases of the stochastic block model and then move on to give practical examples of these ideas on real data sets.
\subsection{The Stochastic Block Model.}
Corollary \ref{corr::extreme_vertices_extremize_attachment} implies that vertices corresponding to the largest-magnitude Fiedler vector entries are also more connected to their own group than the other group, with high probability. Recall the vector
\begin{align*} \Delta = (\# \text{in-cluster neighbors} &- \# \text{out-cluster neighbors})\\
&- \mathbb{E}[\# \text{in-cluster neighbors} - \# \text{out-cluster neighbors}],
\end{align*} 
and recall that it is well-approximated by a Gaussian of standard deviation $\sim_{p,q} \sqrt{n}$. We consider the vector with ``normalized'' standard deviation
\begin{equation}\label{eq::c_def}
     \frac{\Delta}{\sqrt{n}}= \frac{1}{\sqrt{n}}\left((\text{in-degree}-\text{out-degree}) -\frac{(p-q)n}{2} \right)
\end{equation}

\begin{figure}[h!]
\includegraphics[width=\textwidth]{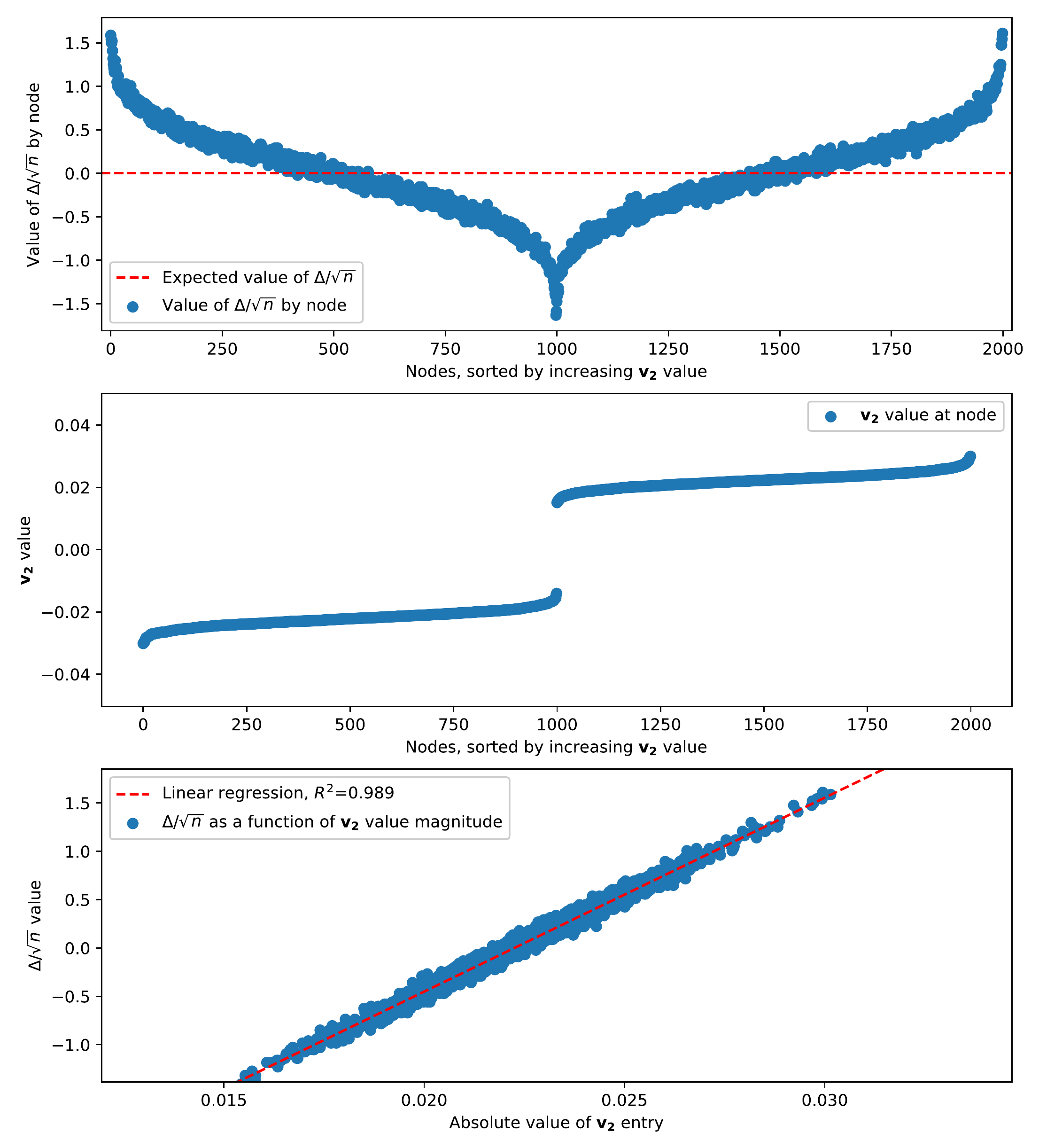}
\caption{(top) The deviation from expected in-group affinity ($c$, defined in Equation \ref{eq::c_def}) for the vertices of a stochastic block model with $(n,p,q) = (2000, 0.6, 0.4)$. Vertices are plotted in increasing order of the corresponding $\mathbf{v_2}$ entry. (mid) Values of $\mathbf{v_2}$ for corresponding vertices, ordered in increasing value. (bottom) Plot showing linear relationship between $\Delta/\sqrt{n}$ and $\left|\mathbf{v_2}\right|$, in accordance with the main theorem.}\label{fig::affinity}
\end{figure}

In Figure \ref{fig::affinity} we plot the values of $\Delta/\sqrt{n}$ against the values of $\mathbf{v_2}$, the eigenvector corresponding to the second largest eigenvalue of the adjacency matrix, for a random instance of the stochastic block model. We order vertices by increasing value of their corresponding $\mathbf{v_2}$ entry, and observe that, as predicted by Corollary \ref{corr::extreme_vertices_extremize_attachment}, the vertices corresponding to the largest-magnitude entries of $\mathbf{v_2}$ have the highest values values of $\Delta/\sqrt{n}$. In particular, with high probability these vertices have a higher in-group affinity than expected, while vertices with a low absolute-value entry have low affinity.
Moreover, our main theorem implies that we would expect a linear relationship between the absolute value of the eigenvector and $\Delta/\sqrt{n}$. Indeed, our numerical experiments empirically validate these findings: the bottom plot on Figure \ref{fig::affinity} displays the results of plotting $\Delta/\sqrt{n}$ against $\left|\mathbf{v_2}\right|$ for an instance of the  stochastic block model.

\subsection{$\mathbf{v_2}$-extreme vertices have a smaller classification error.}
\begin{figure}[h!]
\includegraphics[width=\textwidth]{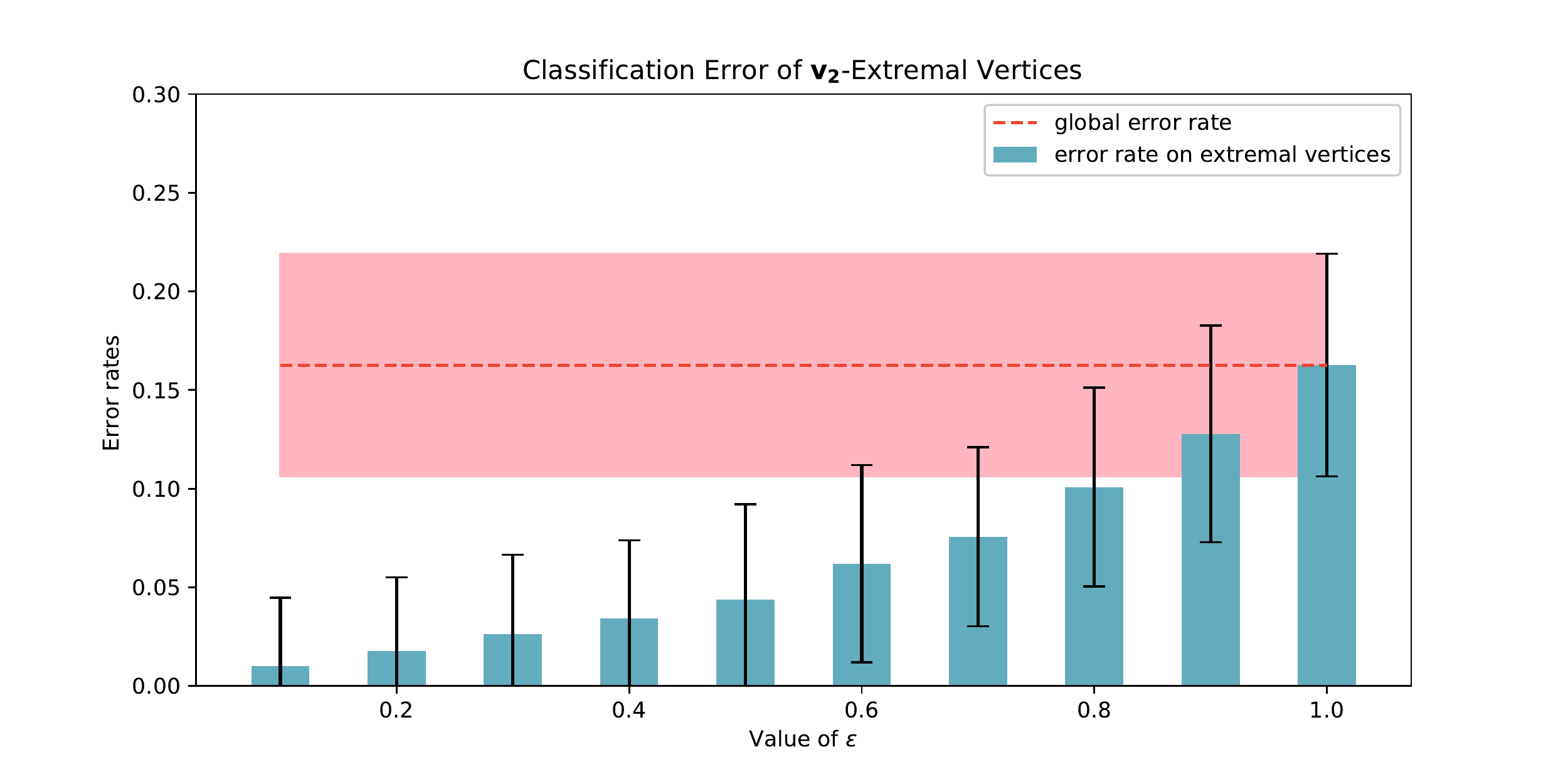}
\caption{Error rates on subsets of vertices with extremal $\mathbf{v_2}$ value, compared with the global $\mathbf{v_2}$ label-estimation error rate. Subsets were chosen by taking the nodes with $\varepsilon\cdot n$ largest magnitude $\mathbf{v_2}$ entries, as in Corollary \ref{corr::extreme_vertices_extremize_attachment}. This figure was generated by randomly sampling 500 independent stochastic block models, $n=200$, $p=0.55,$ and $q=0.45$. }\label{fig::subset_error}
\end{figure}

This observation, that the vertices on which $\mathbf{v_2}$ adopts its extremal values are particularly ``deep'' within their communities, motivates interest in the extremal vertices. In Figure \ref{fig::subset_error}, perform spectral clustering using the signs of $\mathbf{v_2}$, and compare the global error rate to the error rate on the nodes corresponding to the $\varepsilon\cdot n$ largest-magnitude entries of $\mathbf{v_2}$, as in Corollary \ref{corr::extreme_vertices_extremize_attachment}. We demonstrate that with high probability, estimating community identity by the sign of $\mathbf{v_2}$ achieves a lower error rate on sets of extremal vertices.  Additionally,  the error rate on extremal-magnitude vertices is consistently lower than the global error even in ambiguous regimes where the probabilities of in- vs out-community edges are very close. In particular, this gap grows more pronounced as the ratio of the size of the extremal set to the total number of vertices decreases, as illustrated in Figure \ref{fig::subset_error}.

\subsection{MNIST} This consistent outperformance of the $\mathbf{v_2}$ sign estimation method on extremal vertices versus global labeling motivates interest in the usefulness of extremal vertex sets in applications of spectral clustering. Consider the task of distinguishing between classes of visually-similar but symbollically-distinct images: for example, in the classical MNIST dataset, consisting of 28$\times$28 greyscale pixel images of handwritten digits, consider the task of separating the $3$'s and $8$'s into distinct clusters. This task can be formulated as a spectral clustering problem by representing the digits as vertices, and choosing some similarity metric--for example, Euclidean distance--to determine edges between vertices. We ran experiments on MNIST data by adding edges between the vertices corresponding to data points $i$ and $j$ if $i$ is one of the Euclidean-distance $k$-nearest neighbors of $j$, for some specified $k$. In order to ensure symmetry, we make these edges undirected: i.e. if $i$ is one of the $k$-nearest neighbors of $j$, we add an edge from $i$ to $j$ \textit{and} from $j$ to $i$.

\begin{figure}[h!]
\includegraphics[width=\textwidth]{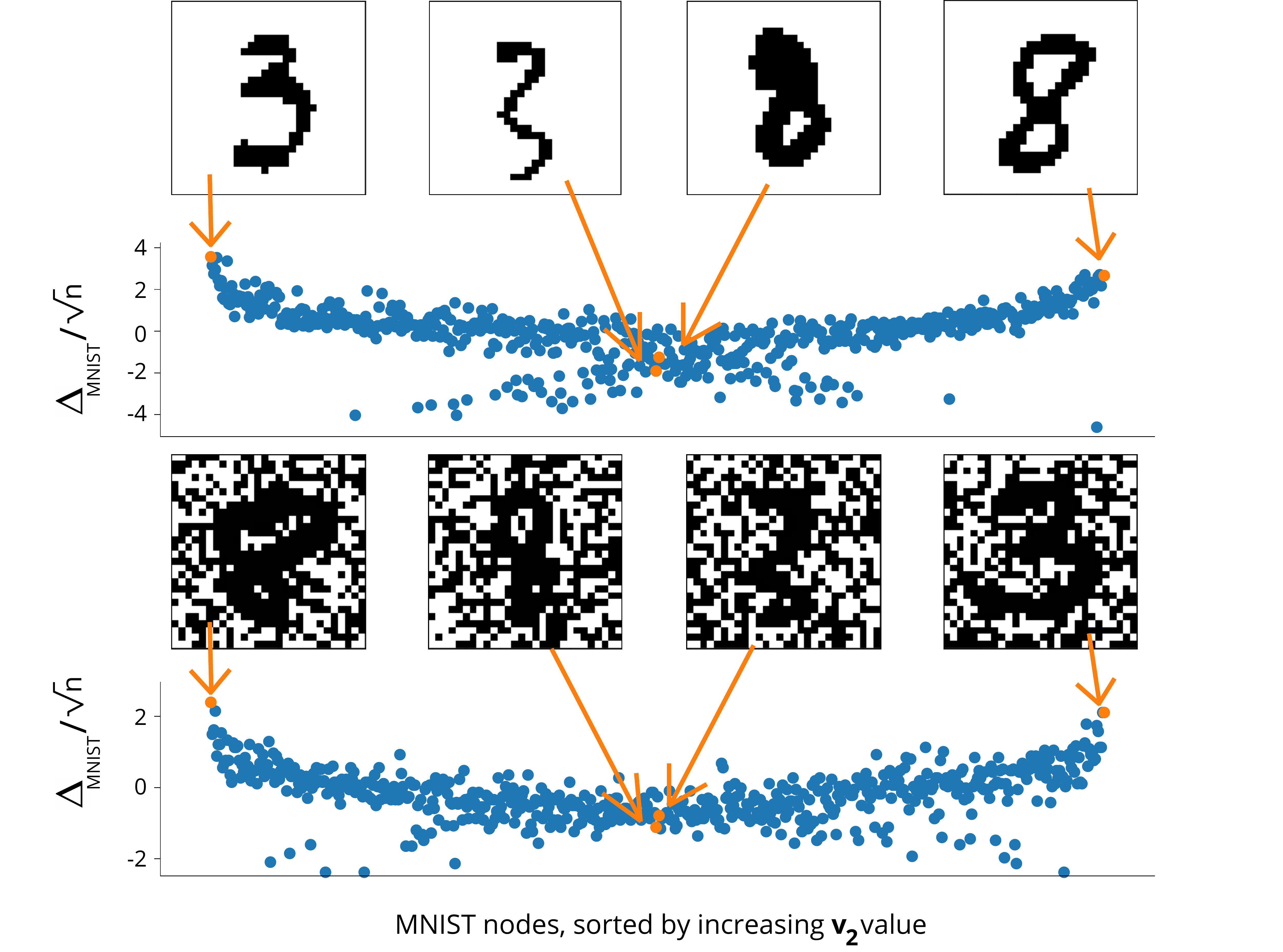}
\caption{Visualization of clustering experiments performed using MNIST dataset. Three hundred images of 3's and three hundred images of 8's were chosen at random from the original MNIST dataset. Pixel values were normalized and rounded to take binary values. A graph was constructed, with a vertex corresponding to each image, and an edge between two vertices if one of the vertices was within the 10\% nearest neighbors of the other, using Euclidean distance. The vector $\mathbf{v_2}$ and values of $c$ (see Equation \ref{eq::c_MNIST}) were calculated for each vertex. The top figure was generated without noise. In the bottom figure, each pixel's binary value was reversed with independent probability $\rho = 0.5$, and the same calculations were performed.}\label{fig::MNIST}
\end{figure}

\indent We observe that in this application, the vector $\mathbf{v_2}$ does indeed encode additional information about the problem. Plotting
\begin{equation}\label{eq::c_MNIST}
    \frac{\Delta_{\text{MNIST}}}{\sqrt{n}} \equiv \frac{1}{\sqrt{n}}\bigg((\# \text{in-cluster neighbors} - \# \text{out-cluster neighbors})-\mu\bigg)
\end{equation}
where $\mu$ is the empirical mean $\# \text{in-cluster neighbors} - \# \text{out-cluster neighbors}$, demonstrates that vertices with extremal $\mathbf{v_2}$values have stronger in-community attachment, as shown in Figure \ref{fig::MNIST}. This agrees with the predictions of Corollary \ref{corr::extreme_vertices_extremize_attachment}. In terms of this application, this trend corresponds to the qualitative observation that 
extremal vertices are associated with more archetypical datapoints, whereas intermediate values correspond to more ambiguous digits, illustrated by visualizing a few extremal and intermediate datapoints in the top panel of Figure \ref{fig::MNIST}. Additionally, this phenomenon is robust under noise: we added noise to the MNIST digits, constructed the nearest-neighbors graph using the same method, and  performed the same calculations of $c_{\text{MNIST}}$. The bottom panel of Figure \ref{fig::MNIST} demonstrates that, even with the presence of noise, vertices with extremal $\mathbf{v_2}$ values have stronger in-community attachment, and are qualitatively more ``easily identifiable'' than vertices with Fiedler value closer to zero. 

\subsection{An Algorithmic Challenge.}
These observations inspire the algorithmic challenge discussed in Section \ref{ssec::main_idea}: Is it possible to utilize the fact that the sign-based spectral clustering error rate for ``extremal'' $\mathbf{v_2}$ entries is much smaller to improve upon the overall classification error? Observing the robustness of this phenonemon even under very high levels of corruption makes this algorithmic challenge even more appealing -- is it possible to propagate the high-quality information (low error in the extremal vertices) to the rest of the vertices?

\section{Proof of the Theorem}
\subsection{Preliminaries}\label{ssec::theorem_prelims}
We recall some basic facts following closely the exposition of Vershynin \cite[\S 4.5]{vershynin}. Let us consider a random adjacency matrix $A$ as induced by the stochastic block model $G(n,p,q)$. It is useful to split these matrices into a deterministic and a random component
$$ A = D + R,~\mbox{where} ~ D = \mathbb{E} A.$$
Assuming the vertices to be nicely ordered so that the first $n/2$ vertices are in the first cluster and the other $n/2$ in the other, this matrix $D$ has a nice form
$$ D = \begin{pmatrix} 
p & p & \dots & p & q &\dots & q & q\\
 p & p & \dots & p & q& \dots & q & q\\
  \dots  & \dots & \dots & \dots & \dots &\dots & \dots & \dots \\
 q & q & \dots & q & p &\dots & p & p\\
  q & q & \dots & q & p& \dots & p & p\\
\end{pmatrix}$$
since it decouples into 4 constant matrices of size $n/2 \times n/2$. (In practice, vertices are not ordered and one tries to recover the order from the entries of the eigenvectors; however, our entire subsequent analysis is invariant under permutation of the entries.) In particular, the matrix $D$ has rank 2, the two eigenvalues are 
$$ \lambda_1(D) = \frac{p+q}{2}n \qquad \mbox{and} \qquad \lambda_2(D) = \frac{p-q}{2}n$$
with corresponding eigenvectors $w_1 = (1, 1, \dots, 1)$ (the constant vector) and $w_2 = (1,1, \dots, 1, -1, \dots, -1, -1)$
(the vector identifying the cluster label). We will henceforth denote these vectors by $\mathbf{1}$ and $\mathbf{g}$. This matrix $D$ is therefore completely understood. The random perturbation $R$ is quite
unpredictable, however, each of its entries is independent of the other and has mean value 0. This shows that (see e.g. \cite{vershynin}) the operator norm
satisfies
$$ \| R \| \leq k \sqrt{n} \qquad \mbox{with likelihood at least}~1 - 4 e^{-n},$$
where $k$ is a universal constant. Recalling the Weyl inequality, we conclude that the $i-$th eigenvalue of $D+R$ can move at most by $\|R\|$: this
means that with high probability, the first two eigenvalues are perturbed by a factor of at most $k\sqrt{n}$ while the remaining eigenvalues are perturbed
away from 0 and are somewhere in the range $[-k\sqrt{n}, k\sqrt{n}]$.
The next ingredient is the Davis-Kahan theorem which we recall for the convience of the reader.
If $M, N \in \mathbb{R}^{n \times n}$ are two symmetric matrices and if the $i-$th eigenvalue of $M$ is well-separated from the rest,
$$ \min_{j \neq i} \left| \lambda_i(M) - \lambda_j(M) \right| = \delta > 0,$$
then the inner product of the $i-$th eigenvalue of $M$ and the $i-$th eigenvalue of $M+N$ satisfies
$$ \sin\left( \angle \mathbf{v_i}(M), \mathbf{v_i}(M+N) \right) \leq \frac{2 \|N\|}{\delta}.$$
We can apply this inequality to $D$ and $R$.  We know that $D$ has an isolated second eigenvalue (the same is true for the first eigenvalue which is also separated), this shows that
$$ \sin\left( \angle \mathbf{v_2}(D), \mathbf{v_2}(D+R) \right) \leq \frac{2 \|R\|}{\delta} \lesssim_{p,q} \frac{1}{\sqrt{n}}.$$
This then implies that most entries of $v_2$ have to have the same sign as the corresponding entry in $v_2(D)$ which we know to be the ground truth and
this is how we see the validity of spectral clustering. In our subsequent argument, we will also need some information on the first eigenvector $\mathbf{v_1}$, by the same argument
$$ \sin\left( \angle \mathbf{v_1}(D), \mathbf{v_1}(D+R) \right) \leq \frac{2 \|R\|}{\delta} \lesssim_{p,q} \frac{1}{\sqrt{n}}.$$

\subsection{A Lemma.} Before embarking on a proof of the main result, we derive a Lemma that contains the bulk of the argument. Let us again return to considering $A$ as the
adjacency matrix of a matrix from the stochatic block model. The Davis-Kahan argument outlined in \S \ref{ssec::theorem_prelims} allows us to assume that the first
eigenvalue has (up to scaling) the form
$$ \mathbf{v_1} = \frac{\mathbf{1}}{\sqrt{n}} + \mathbf{e_1},~\qquad \mbox{where the error satisfies}~\| \mathbf{e_1}  \| \lesssim \frac{1}{\sqrt{n}} \qquad \mbox{w.h.p.},$$
where $\mathbf{e_1}$ is orthogonal to $\mathbf{1}$ and $\mathbf{1} = (1,1,\dots, 1)$. Moreover, by the same reasoning, introducing the group identification vector $\textbf{g} = (1,1,\dots, 1, -1, \dots, -1, -1)$,
we get that
$$ \mathbf{v_2} = \frac{\mathbf{g}}{\sqrt{n}} + \mathbf{e_2},~\qquad \mbox{where the error satisfies}~\| \mathbf{e_2}  \| \lesssim \frac{1}{\sqrt{n}} \quad \mbox{w.h.p.}$$
where $\mathbf{e_2}$ is orthogonal to $\mathbf{g}$.
We have
$$ \left\langle \mathbf{v_1}, \mathbf{e_1} \right\rangle = \left\langle  \frac{\mathbf{1}}{\sqrt{n}} + \mathbf{e_1}, \mathbf{e_1} \right\rangle = \|\mathbf{e_1}\|^2 \lesssim \frac{1}{n} \qquad \mbox{w.h.p.}$$
and, by the same reasoning,
$$ \left\langle \mathbf{v_2}, \mathbf{e_2} \right\rangle = \left\langle  \frac{\mathbf{g}}{\sqrt{n}} + \mathbf{e_2}, \mathbf{e_2} \right\rangle = \|\mathbf{e_2}\|^2 \lesssim \frac{1}{n} \qquad \mbox{w.h.p.}$$

We will now establish the crucial technical ingredient. We will show that the perturbation $\mathbf{e_2}$ added to the leading term $\mathbf{g}/\sqrt{n}$ for the second eigenvector
is almost orthogonal to the leading eigenvector $\mathbf{v_1}$. The proof is overall quite simple: we decompose all ingredients into their basic
building blocks and, having made a good ansatz, most of the arising computations turn out to be easy to deal with. There are two large terms which reduce
to simple properties of the random graph that can be completely analyzed using the central limit theorem.

\begin{lemma} For every $\varepsilon_0>0$, for $n$ sufficiently large, we have
$$ \left| \left\langle \mathbf{v_1}, \mathbf{e_2} \right\rangle \right| \leq n^{-1 + \varepsilon_0}  \qquad \mbox{with high probability.}$$
\end{lemma}
The proof allows for slightly refined estimates. In practice, we show that the quantity behaves as $\sim n^{-1}\cdot \mathcal{N}(0,1)$, where $\mathcal{N}(0,1)$ is a standard Gaussian. This is clearly less than $n^{-1 + \varepsilon_0}$ with high probability. 

\begin{proof} Let us suppose that  
$ \left| \left\langle \mathbf{v_1}, \mathbf{e_2} \right\rangle \right|$ is large.
We observe that
\begin{align*}
\left| \left\langle  \mathbf{v_1}, \mathbf{e_2} \right\rangle \right|&= \left| \left\langle \mathbf{v_1}, \frac{\mathbf{g}}{\sqrt{n}} + \mathbf{e_2} \right\rangle -  \left\langle \mathbf{v_1}, \frac{\mathbf{g}}{\sqrt{n}} \right\rangle \right|\\
&=  \left| \left\langle  \mathbf{v_1}, \mathbf{v_2} \right\rangle -  \left\langle \frac{\mathbf{1}}{\sqrt{n}} + \mathbf{e_1}, \frac{\mathbf{g}}{\sqrt{n}} \right\rangle \right| = \left| \left\langle \mathbf{e_1}, \frac{\mathbf{g}}{\sqrt{n}}\right\rangle \right|.
\end{align*}
If $ \left| \left\langle \mathbf{v_1}, \mathbf{e_2} \right\rangle \right|$ were large, it would have the interesting implication that the small perturbation $\mathbf{e_1}$ on top of our guess for the leading eigenvector, $\mathbf{1}$, is strongly correlated with the leading term for the second eigenvector. The rest of the proof is concerned with showing that this is not the case.
We rewrite $\mathbf{v_1}$ in a way that makes the intrinsic dependence of the error term in terms of the second eigenvector clearer and  define $\mathbf{r}\in \mathbb{R}^n$ as the vector satisfying the equation
$$ \mathbf{ v_1} = \frac{\mathbf{1}}{\sqrt{n}} + \frac{c}{\sqrt{n}} \frac{\mathbf{g}}{\sqrt{n}} + \mathbf{r},$$
where the remainder $\mathbf{r}$ is orthogonal to $\mathbf{1}$ and $\mathbf{g}$ (this, implicitly, also defines the value of $c$ by orthogonality). By the Davis-Kahan error bound,
\[
    \norm{\frac{c}{\sqrt{n}} \frac{\mathbf{g}}{\sqrt{n}} + \mathbf{r}} \lesssim n^{-1/2} \qquad \text{w.h.p.}
\]
and $\mathbf{r}$ and $\mathbf{g}$ are orthogonal, so by Pythagorean theorem
\[
    \norm{\mathbf{r}}^2 \leq \norm{\frac{c}{\sqrt{n}} \frac{\mathbf{g}}{\sqrt{n}} + \mathbf{r}}^2 \lesssim n^{-1},
\]
so $\norm{\mathbf{r}} \lesssim n^{-1/2}$.
We want to show that $ \left| \left\langle \mathbf{v_1}, \mathbf{e_2} \right\rangle \right| $ is small and, by the argument above, it suffices to show that $| \left\langle \mathbf{e_1}, \mathbf{g}/\sqrt{n} \right\rangle|$ is small. Observe that 
$$ \left\langle \mathbf{e_1}, \frac{\mathbf{g}}{\sqrt{n}}\right\rangle =  \left\langle \frac{c}{\sqrt{n}} \frac{\mathbf{g}}{\sqrt{n}} + \mathbf{r}, \frac{\mathbf{g}}{\sqrt{n}}\right\rangle  =  \frac{c}{\sqrt{n}}.$$
Our goal is now to show that $|c|$ has to be small. By definition, $\mathbf{v_1}$ attains the largest Rayleigh quotient
$$\lambda_1 =  \frac{\left\langle \mathbf{v_1}, A \mathbf{v_1} \right\rangle}{ \left\langle \mathbf{v_1}, \mathbf{v_1} \right\rangle} = \sup_{\mathbf{w} \neq 0}   \frac{\left\langle \mathbf{w}, A \mathbf{w} \right\rangle}{ \left\langle \mathbf{w}, \mathbf{w} \right\rangle} $$
among all vectors of a fixed size.
We compare the size of this Rayleigh quotient to the Rayleigh quotient for the `competing' vector
$$ \mathbf{w} = \sqrt{1+\frac{c^2}{n}}  \frac{\mathbf{1}}{\sqrt{n}}   +\ \mathbf{r}.$$
 The motivation behind this construction of $\mathbf{w}$ is quite simple: $\mathbf{1}/\sqrt{n}$ is the leading term for $\mathbf{v_1}$, whereas $\mathbf{g}/\sqrt{n}$ is the leading term for $\mathbf{v_2}$, which has a smaller Rayleigh quotient. It thus seems reasonable to wonder whether we can increase the Rayleigh quotient by re-distributing components of $\mathbf{v_1}$ from $\mathbf{g}$ to $\mathbf{1}$. Coefficients are then chosen to ensure $\|\mathbf{w}\| = \|\mathbf{v_1}\|$. We compute 
\begin{align*}
\left\langle \mathbf{w}, A \mathbf{w} \right\rangle &= \left(1+\frac{c^2}{n}\right)  \left\langle  \frac{\mathbf{1}}{\sqrt{n}} , A  \frac{\mathbf{1}}{\sqrt{n}}\right\rangle +   \left\langle \sqrt{1+\frac{c^2}{n}}  \frac{\mathbf{1}}{\sqrt{n}} , A \mathbf{r}\right\rangle\\
& +\left\langle\ \mathbf{r}, A\sqrt{1+\frac{c^2}{n}}  \frac{\mathbf{1}}{\sqrt{n}}  \right\rangle +  \left\langle\ \mathbf{r}, A\mathbf{r}  \right\rangle.
\end{align*}
We recall that $A$ is symmetric allowing us to combine two terms and write
\begin{align*}
\left\langle \mathbf{w}, A \mathbf{w} \right\rangle &= \left(1+\frac{c^2}{n}\right)  \left\langle  \frac{\mathbf{1}}{\sqrt{n}} , A  \frac{\mathbf{1}}{\sqrt{n}}\right\rangle +  2 \left\langle \sqrt{1+\frac{c^2}{n}}  \frac{\mathbf{1}}{\sqrt{n}} , A \mathbf{r}\right\rangle  +  \left\langle\ \mathbf{r}, A\mathbf{r}  \right\rangle.
\end{align*}

Likewise, we have
\begin{align*}
\left\langle \mathbf{v_1}, A \mathbf{v_1} \right\rangle &=  \left\langle  \frac{\mathbf{1}}{\sqrt{n}} , A  \frac{\mathbf{1}}{\sqrt{n}}\right\rangle + 2\left\langle \frac{\mathbf{1}}{\sqrt{n}}, A \frac{c}{\sqrt{n}} \frac{\mathbf{g}}{\sqrt{n}} \right\rangle\\
&+ 2\left\langle \frac{\mathbf{1}}{\sqrt{n}}, A \mathbf{r} \right\rangle + \left\langle   \frac{c}{\sqrt{n}} \frac{\mathbf{g}}{\sqrt{n}}, A \frac{c}{\sqrt{n}} \frac{\mathbf{g}}{\sqrt{n}} \right\rangle \\
&  +2\left\langle   \frac{c}{\sqrt{n}} \frac{\mathbf{g}}{\sqrt{n}}, A \mathbf{r} \right\rangle
+  \left\langle\ \mathbf{r}, A\mathbf{r}  \right\rangle.
\end{align*}
Since $\mathbf{v_1}$ is the largest eigenvector, we have
$$ Q \equiv \left\langle \mathbf{w}, A \mathbf{w} \right\rangle - \left\langle \mathbf{v_1}, A \mathbf{v_1} \right\rangle  \leq 0.$$
We will derive a lower bound on $Q$ depending on $c$: this lower bound will then imply that $c$ has to be small. We start by writing out $Q$.
\begin{align*}
Q &= \frac{c^2}{n}  \left\langle  \frac{\mathbf{1}}{\sqrt{n}} , A  \frac{\mathbf{1}}{\sqrt{n}}\right\rangle + 2\left( \sqrt{1 + \frac{c^2}{n}} - 1 \right)\left\langle \frac{\mathbf{1}}{\sqrt{n}}, A \mathbf{r} \right\rangle\\
&- 2\left\langle \frac{\mathbf{1}}{\sqrt{n}}, A \frac{c}{\sqrt{n}} \frac{\mathbf{g}}{\sqrt{n}} \right\rangle -  \left\langle   \frac{c}{\sqrt{n}} \frac{\mathbf{g}}{\sqrt{n}}, A \frac{c}{\sqrt{n}} \frac{\mathbf{g}}{\sqrt{n}} \right\rangle\ -2\left\langle   \frac{c}{\sqrt{n}} \frac{\mathbf{g}}{\sqrt{n}}, A \mathbf{r} \right\rangle.
\end{align*}
We will discuss this quantity $Q$ term-by-term. \\

\textit{The first term.} We observe that 
$$ \frac{c^2}{n}  \left\langle  \frac{\mathbf{1}}{\sqrt{n}} , A  \frac{\mathbf{1}}{\sqrt{n}}\right\rangle = \frac{c^2}{n^2} \cdot \sum_{i=1}^n \text{degree}(i) = \frac{c^2}{n^2} \cdot 2\cdot (\# \mbox{edges in the graph}).$$
By the central limit theorem, this number is tightly concentrated and 
$$ \frac{c^2}{n^2} \cdot 2 \cdot (\# \mbox{edges in the graph}) = c^2(p+q) + \mathcal{O}\left(\frac{1}{n}\right) \qquad \mbox{w.h.p.}$$

\textit{The second term.} The second term is small. Taylor expansion shows that
$$ 2\left( \sqrt{1 + \frac{c^2}{n}} - 1 \right) \sim \frac{c^2}{n} + \mbox{l.o.t.}$$
As for the inner product, by Cauchy Schwarz
$$ \left| \left\langle \frac{\mathbf{1}}{\sqrt{n}}, A \mathbf{r} \right\rangle\right| \leq \left\|\frac{\mathbf{1}}{\sqrt{n}} \right\| \| A \mathbf{r}\| = \| A \mathbf{r}\|.$$
Using the spectral theorem, we obtain
$$ \| A \mathbf{r}\|^2 = \sum_{k=1}^{n} \lambda_k^2 \left\langle \mathbf{r} , \mathbf{v_k} \right\rangle^2.$$
We know that the first two eigenvalues are large, $\lambda_1 \sim n \sim \lambda_2$, the remaining eigenvalues are all at scale $\sim \sqrt{n}$. As a consequence of orthogonality,
$$ \left\langle \mathbf{r} , \mathbf{v_1} \right\rangle^2 = \left\langle \mathbf{r} , \mathbf{r} \right\rangle^2 = \|\mathbf{r}\|^4 \lesssim \frac{1}{n^2}.$$
Likewise, we have
$$ \left\langle \mathbf{r} , \mathbf{v_2} \right\rangle^2 = \left\langle \mathbf{r} , \frac{\mathbf{g}}{\sqrt{n}} + \mathbf{e_2} \right\rangle^2 =  \left\langle \mathbf{r} , \mathbf{e_2} \right\rangle^2 \leq \|\mathbf{r}\|^2 \| \mathbf{e_2}\|^2 \lesssim \frac{1}{n^2}.$$
As for the remaining terms, we use the Pythagorean theorem to write
$$ \sum_{k=3}^{n} \lambda_k^2 \left\langle \mathbf{r} , \mathbf{v_k} \right\rangle^2 \leq \left(\max_{3 \leq k \leq n}{ \lambda_k^2}\right)  \sum_{k=3}^{n} \left\langle \mathbf{r} , \mathbf{v_k} \right\rangle^2 \lesssim n \|\mathbf{r}\|^2 \lesssim 1.$$
Altogether, we get that the second term can be bounded by
$$ \left| \left( \sqrt{1 + \frac{c^2}{n}} - 1 \right)\left\langle \frac{\mathbf{1}}{\sqrt{n}}, A \mathbf{r} \right\rangle \right|
\lesssim \frac{c^2}{n} \| A \mathbf{r}\| \lesssim \frac{c^2}{n}.$$

\textit{The third term.} We write the term, using the symmetry of $A$, as
\begin{align*}
\left\langle \frac{\mathbf{1}}{\sqrt{n}}, A \frac{c}{\sqrt{n}} \frac{\mathbf{g}}{\sqrt{n}} \right\rangle =  \frac{c}{n^{3/2}} \left\langle \mathbf{1}, A \mathbf{g} \right\rangle  = \frac{c}{n^{3/2}} \left\langle A \mathbf{1},  \mathbf{g} \right\rangle.
\end{align*}
However, $A \mathbf{1}$ merely counts the degree of each vertex. In total, the inner product therefore computes the sum of all degree of vertices in the first cluster and subtracts the sum of the degrees of the vertices in the second cluster. Both numbers are the same in expectation. Moreover, both sums are asymptotically distributed like a Gaussian with expectation
\begin{align*} 
    \mathbb{E} ~\sum_{i=1}^{n/2} \sum_{j=1}^{n}{ a_{ij}}  &= \sum_{i=1}^{n/2} \left( \sum_{j=1}^{n/2}{ \mathbb{E}~a_{ij}}+\sum_{j=n/2+1}^{n}{\mathbb{E}~a_{ij}}\right)\\
    &=  \sum_{i=1}^{n/2} \left( \sum_{j=1}^{n/2}{ p}+\sum_{j=n/2+1}^{n}{q}\right) = \left(\frac{n}{2}\right)^2 (p+q).
\end{align*}
When we subtract them, the expectation is 0 and it remains to control for the variance, which is easily seen to be
$$
 \mathbb{V} \sum_{i=1}^{n/2} \sum_{j=1}^{n}{ a_{ij}} \sim_{p,q} n^2.$$
Accounting for the rescaling by the $n^{-3/2}$ factor, the third term behaves like a Gaussian centered at 0 having standard deviation $n^{-1}$. In particular, we can expect for $n$ large that
$$  \frac{c}{n^{3/2}} \left\langle A \mathbf{1},  \mathbf{g} \right\rangle \sim_{p,q} \frac{c}{\sqrt{n}} \mathcal{N}(0,1),$$
where the last $\sim$ is to be understand in the sense of indicating scale. The quantity might be positive or negative depending on the particular instance but the bulk of the probability is accurate modeled by this rescaled Gaussian.\\

\textit{The fourth term.} The fourth term is not small. Indeed, we can write
\begin{align*}
 A \frac{c}{\sqrt{n}} \frac{\mathbf{g}}{\sqrt{n}} &= A \left( \frac{c}{\sqrt{n}} \frac{\mathbf{g}}{\sqrt{n}} + \frac{c}{\sqrt{n}} \mathbf{e_2} \right) - A  \frac{c}{\sqrt{n}} \mathbf{e_2} \\
&=A  \frac{c}{\sqrt{n}} \mathbf{v_2} - A  \frac{c}{\sqrt{n}} \mathbf{e_2}  = \lambda_2 \frac{c}{\sqrt{n}} \mathbf{v_2} - A  \frac{c}{\sqrt{n}} \mathbf{e_2}.
 \end{align*}
 This shows that
 \begin{align*}
  \left\langle   \frac{c}{\sqrt{n}} \frac{\mathbf{g}}{\sqrt{n}}, A \frac{c}{\sqrt{n}} \frac{\mathbf{g}}{\sqrt{n}} \right\rangle = \lambda_2  \left\langle   \frac{c}{\sqrt{n}} \frac{\mathbf{g}}{\sqrt{n}},  \frac{c}{\sqrt{n}}\mathbf{v_2} \right\rangle -  \left\langle   \frac{c}{\sqrt{n}} \frac{\mathbf{g}}{\sqrt{n}},A  \frac{c}{\sqrt{n}} \mathbf{e_2}\right\rangle.
 \end{align*}
 Here, the first term is large since, by orthogonality of $\mathbf{g}$ and $\mathbf{e_2}$,
 \begin{align*}
 \lambda_2  \left\langle   \frac{c}{\sqrt{n}} \frac{\mathbf{g}}{\sqrt{n}},  \frac{c}{\sqrt{n}}\mathbf{v_2} \right\rangle =  \lambda_2  \left\langle   \frac{c}{\sqrt{n}} \frac{\mathbf{g}}{\sqrt{n}},   \frac{c}{\sqrt{n}}\frac{\mathbf{g}}{\sqrt{n}} \right\rangle = \lambda_2 \frac{c^2}{n}.
 \end{align*}
 We observe that we have explicit bounds on $\lambda_2$ and
 $$ \lambda_2 \frac{c^2}{n} = \left( \frac{p-q}{2} n + \mathcal{O}(\sqrt{n})\right) \frac{c^2}{n} = c^2 \left(\frac{p-q}{2}\right) + \mathcal{O}\left(\frac{c^2}{\sqrt{n}}\right) \quad \mbox{w.h.p.}$$
It remains to show that the second term is small. We note that 
\begin{align*}
  \left\langle   \frac{c}{\sqrt{n}} \frac{\mathbf{g}}{\sqrt{n}},A  \frac{c}{\sqrt{n}} \mathbf{e_2}\right\rangle &= \frac{c^2}{n}   \left\langle   A \frac{\mathbf{g}}{\sqrt{n}},   \mathbf{e_2}\right\rangle \\
  &= \frac{c^2}{n}   \left\langle   A\left( \frac{\mathbf{g}}{\sqrt{n}} + \mathbf{e_2} \right),   \mathbf{e_2}\right\rangle - \frac{c^2}{n}   \left\langle   A \mathbf{e_2},  \mathbf{e_2}\right\rangle \\
    &= \frac{c^2}{n}   \left\langle   A\mathbf{v_2},   \mathbf{e_2}\right\rangle - \frac{c^2}{n}   \left\langle   A \mathbf{e_2},   \mathbf{e_2}\right\rangle\\
    &=    \frac{c^2}{n} \lambda_2  \left\langle   \mathbf{v_2} ,   \mathbf{e_2}\right\rangle - \frac{c^2}{n}   \left\langle   A \mathbf{e_2},   \mathbf{e_2}\right\rangle \\
    &= \frac{c^2}{n} \lambda_2   \|\mathbf{e_2}\|^2 - \frac{c^2}{n}   \left\langle   A \mathbf{e_2},   \mathbf{e_2}\right\rangle.
\end{align*}
We recall that $\|\mathbf{e_2}\| \lesssim n^{-1/2}$ and $\lambda_2 \sim_{p,q} n$ thus
$$ \frac{c^2}{n} \lambda_2   \|\mathbf{e_2}\|^2 \lesssim_{p,q} \frac{c^2}{n}.$$
As for the second term, we have the trivial estimate
$$ \left| \frac{c^2}{n}   \left\langle   A \mathbf{e_2},   \mathbf{e_2}\right\rangle \right| \leq \frac{c^2}{n} \|A\|   \|\mathbf{e_2}\|^2 = \frac{c^2}{n}\lambda_1  \|\mathbf{e_2}\|^2 \lesssim  \frac{c^2}{n}.$$

\textit{The fifth term.} It remains to estimate the fifth term, which is small. We use
\begin{align*}
 \left\langle   \frac{c}{\sqrt{n}} \frac{\mathbf{g}}{\sqrt{n}}, A \mathbf{r} \right\rangle &= \frac{c}{\sqrt{n}} \left\langle A\frac{\mathbf{g}}{\sqrt{n}},  \mathbf{r} \right\rangle \\
 &= \frac{c}{\sqrt{n}} \left\langle A\left( \frac{\mathbf{g}}{\sqrt{n}} + \mathbf{e_2} \right),  \mathbf{r} \right\rangle -  \frac{c}{\sqrt{n}} \left\langle A \mathbf{e_2},  \mathbf{r} \right\rangle \\
  &= \frac{c}{\sqrt{n}} \left\langle A \mathbf{v_2},  \mathbf{r} \right\rangle -  \frac{c}{\sqrt{n}} \left\langle A \mathbf{e_2},  \mathbf{r} \right\rangle \\
  &= \frac{c}{\sqrt{n}} \lambda_2 \left\langle \mathbf{v_2},  \mathbf{r} \right\rangle -  \frac{c}{\sqrt{n}} \left\langle  \mathbf{e_2},  A\mathbf{r} \right\rangle \\
    &= \frac{c}{\sqrt{n}} \lambda_2 \left\langle \mathbf{e_2},  \mathbf{r} \right\rangle -  \frac{c}{\sqrt{n}} \left\langle  \mathbf{e_2},  A\mathbf{r} \right\rangle.
 \end{align*}
 The first term is small since 
 \begin{align*}
  \frac{c}{\sqrt{n}} \lambda_2 \left| \left\langle \mathbf{e_2},  \mathbf{r} \right\rangle \right| \lesssim \frac{c}{\sqrt{n}} n \| \mathbf{e_2}\| \| \mathbf{r}\| \lesssim \frac{c}{\sqrt{n}}.
 \end{align*}
 It remains to understand the term $\left\langle \mathbf{e_2}, A\mathbf{r} \right\rangle$. To this end, we argue that 
 \begin{align*}
    \|A\mathbf{r}\|^2 &= \lambda_1(A)^2 \left\langle \mathbf{r}, \mathbf{v_1} \right\rangle^2 + \lambda_2(A)^2 \left\langle \mathbf{r}, \mathbf{v_2} \right\rangle^2 
    +\sum_{i\geq3}^n\lambda_i(A)^2 \langle \mathbf{r},\mathbf{v_i}\rangle^2  \\
    &\leq \lambda_1(A)^2 \left\langle \mathbf{r}, \mathbf{v_1} \right\rangle^2 + \lambda_2(A)^2 \left\langle \mathbf{r}, \mathbf{v_2} \right\rangle^2 + \lambda_3(A)^2\sum_{i=3}^{n}\langle \mathbf{r},\mathbf{v_i}\rangle^2\\
    &\leq \lambda_1(A)^2 \left\langle \mathbf{r}, \mathbf{v_1} \right\rangle^2 + \lambda_2(A)^2 \left\langle \mathbf{r}, \mathbf{v_2} \right\rangle^2 + \lambda_3(A)^2\norm{\mathbf{r}}^2
 \end{align*}
 
 However, by orthogonality,
 $$ \left\langle \mathbf{r}, \mathbf{v_1} \right\rangle^2 = \left\langle \mathbf{r}, \mathbf{r} \right\rangle^2 = \|\mathbf{r}\|^4 \lesssim \frac{1}{n^2}.$$
As for the second term, we recall that $\mathbf{r}$ is orthogonal to $\mathbf{g}$ and therefore
$$  \left\langle \mathbf{r}, \mathbf{v_2} \right\rangle^2 =  \left\langle \mathbf{r}, \mathbf{e_2} \right\rangle^2 \leq \| \mathbf{r}\|^2 \|\mathbf{e}_2\|^2 \lesssim \frac{1}{n^2}.$$
Moreover, we have $\lambda_3(A) \lesssim \sqrt{n}$ and thus the third term satisfies
$$ \lambda_3(A)^2 \norm{\mathbf{r}}^2 \lesssim 1.$$
 Altogether, we obtain 
$$ \left| \frac{c}{\sqrt{n}} \left\langle  \mathbf{e_2},  A\mathbf{r} \right\rangle \right| \leq \frac{c}{\sqrt{n}} \| \mathbf{e_2}\| \| A \mathbf{r} \| \lesssim \frac{c}{n}.$$

\textit{Summary.} These five estimates come with two different types of guarantees. They are all probabilistic. The second, fourth and fifth term come with explicit constants that hold with high probability (i.e. a probability converging exponentially in $n$ to 1). The first comes with a Gaussian error that acts on a small scale $n^{-1}$. The error term in the third quantity is Gaussian and at scale $\sim n^{-1/2}$. 
This means that we have obtained, by an abuse of notation, the estimate $$ \alpha_{1} c^2 -  \alpha_2 \frac{c^2}{n} - \frac{c}{\sqrt{n}} \mathcal{N}(0,1) - \alpha_3 \frac{c}{n} \leq Q \leq 0$$
where $\alpha_{1}, \alpha_2, \alpha_3 >0$ are fixed constants and $\mathcal{N}(0,1)$ is a standard Gaussian. Thus
\begin{equation}\label{eq::c_distribution}
     c^2  \lesssim_{\alpha_1,\alpha_2,\alpha_3}  \frac{c^2}{n} + \frac{c}{\sqrt{n}} \mathcal{N}(0,1) + \frac{c}{n}.
\end{equation} 
Recalling that, by orthogonality
$$ \left\langle \mathbf{e_1}, \frac{\mathbf{g}}{\sqrt{n}}\right\rangle =  \left\langle \frac{c}{\sqrt{n}} \frac{\mathbf{g}}{\sqrt{n}} + \mathbf{r}, \frac{\mathbf{g}}{\sqrt{n}}\right\rangle  =  \frac{c}{\sqrt{n}},$$
we see that Equation \ref{eq::c_distribution} leads to a contradiction with high probability as soon as $|c| \gtrsim n^{-1/2 + \varepsilon_0}$. We thus conclude that
$$
    \left|\langle \mathbf{v_1},\mathbf{e_2}\rangle \right| = \left|\langle \mathbf{e_1},\frac{\mathbf{g}}{\sqrt{n}}\rangle \right| = \frac{|c|}{\sqrt{n}} \lesssim n^{-1 + \varepsilon_0}.
$$
\end{proof}

\subsection{Proof of the Main Theorem}
\begin{proof} We will expand the eigenvalue equation
\begin{equation} \label{eineq}
A \mathbf{v_2} = \lambda_2 \mathbf{v_2}
 \end{equation}
As above, we make the ansatz
$$ \mathbf{v_2} = \frac{\mathbf{g}}{\sqrt{n}} + \mathbf{e_2}$$
and we know that the error satisfies
$$\| \mathbf{e_2}  \| \lesssim \frac{1}{\sqrt{n}} \quad \mbox{w.h.p.}$$
We also introduce the constant $k \in \mathbb{R}$ by writing
\begin{equation} \label{ansatz2}
 \lambda_2 = \frac{p-q}{2}n + k \sqrt{n}.
 \end{equation}
 We have $|k| \lesssim 1$ with high probability. We introduce $\mathbf{dev} \in \mathbb{R}^n$ by writing
\begin{equation}  \label{ansatz3}
A\mathbf{g}= \frac{(p-q)n}{2}\mathbf{g}+ \mathbf{dev}.
\end{equation}
Equipped with the definitions \eqref{ansatz2} and \eqref{ansatz3}, we can re-expand \eqref{eineq} and obtain
$$ \frac{1}{\sqrt{n}} \mathbf{dev} + A \mathbf{e_2} = k \mathbf{g} + \frac{p-q}{2}n \mathbf{e_2} + k \sqrt{n}\mathbf{e_2}.$$
 We rewrite this equation as
 $$ \left(  \mbox{Id}_{n \times n} +  \frac{2 k}{(p-q)\sqrt{n}}\cdot \mbox{Id}_{n \times n}\right) \mathbf{e_2} =  \frac{2 \cdot \mathbf{dev}}{(p-q)n^{3/2}}   - \frac{2 k \mathbf{g} }{(p-q)n} +  \frac{2A }{(p-q)n}\mathbf{e_2}$$
We show that the term $A\mathbf{e_2}$ on the right-hand side is going to be small. We use the spectral theorem to estimate
$$ \|A \mathbf{e_2}\|^2 \leq \lambda_1^2 \left\langle \mathbf{e_2}, \mathbf{v_1} \right\rangle^2 +  \lambda_2^2 \left\langle \mathbf{e_2}, \mathbf{v_2} \right\rangle^2 + \lambda_3^2 \|\mathbf{e_2}\|^2.$$
Recalling Lemma 1, we get, for every $\varepsilon_0 > 0$ and with high probability for $n$ sufficiently large (depending on $\varepsilon_0$),
$$ \|A \mathbf{e_2}\|^2 \lesssim \frac{  \lambda_1^2}{n^{2-2\varepsilon_0}} + 2\lambda_2^2 \cdot \|\mathbf{e_2}\|^4 + 1 \lesssim n^{2\varepsilon_0}$$
and thus
 $$ \left\| \frac{2A }{(p-q)n}\mathbf{e_2}\right\| \lesssim \frac{1}{n^{1-\varepsilon_0}}.$$
However, this is quite small and much smaller than the other terms on the right-hand side. Using a Neumann series (see e.g. Kato \cite{kato})
 $$ \left( \mbox{Id}_{n \times n} - T\right)^{-1} = \mbox{Id}_{n \times n} + \sum_{\ell=1}^{\infty} T^\ell$$
 in combination with the basic estimate
 $$ \left\|  \sum_{\ell=1}^{\infty} T^\ell \right\| \leq  \sum_{\ell=1}^{\infty} \|T\|^\ell = \frac{\|T\|}{1 - \|T\|}$$
 applied to
 $$ T =  \frac{2 k}{(p-q)\sqrt{n}}\cdot \mbox{Id}_{n \times n} \qquad \mbox{satisfying, with high probability,} \qquad \|T\| \lesssim \frac{1}{n},$$
 we can invert the matrix on the left-hand side and obtain
 $$ \mathbf{e_2} = \frac{2\cdot \mathbf{dev}}{(p-q)n^{3/2}}   - \frac{2 k \mathbf{g}}{(p-q)n} + \mathcal{O}\left( \frac{1}{n^{1-\varepsilon_0}}\right).$$
 Recalling the definitions
$$  \mathbf{dev} =A\mathbf{g} - \frac{(p-q)n}{2}\mathbf{g} \qquad \mbox{and} \qquad  k = \frac{\lambda_2 - \frac{p-q}{2}n}{\sqrt{n}}$$
we arrive at the desired statement.
\end{proof}

\subsection{Proof of Corollary 1}
\begin{proof} We recall that
$$ \mathbf{v_2} = (1+\varepsilon_{\text{global}})\frac{\mathbf{g}}{\sqrt{n}}  + \frac{2}{(p-q)n^{3/2}}   \left(A\mathbf{g} - \frac{(p-q)n}{2}\mathbf{g}\right)  + \mathbf{error}$$
where $\mathbf{error}$ is of size $\|\mathbf{error}\| \lesssim n^{-1 + \varepsilon_0}$ and $\varepsilon_{\text{global}}$ corresponds to the global shift, discussed in \S \ref{ssec::global_shift}. 
We will only carry out the argument for the largest entries of $\mathbf{v_2}$: the argument for the smallest entries of $\mathbf{v_2}$ is identical up to changes of sign. Ignoring $\mathbf{error}$ for now, we will try to understand the second term. It measures the difference between in-group connections and out-group connections from their expected size. For a single vertex, this term is approximately gaussian and $\sim c_{p,q} \cdot \sqrt{n} \cdot \mathcal{N}(0,1)$ for some constant $c_{p,q} > 0$. 
Let us now consider the likelihood that this Gaussian random variable is bigger than $2\eta$. This likelihood, $\varepsilon$, is bigger than 0 (though decaying quite rapidly as $\eta$ increasing). This means, that we asymptotically expect $\varepsilon \cdot n$ of the vertices to exceed expectation by $c_{p,q}\cdot \sqrt{n} \cdot (2\eta)$. We will now argue that $99\%$ of the largest $\varepsilon \cdot n$ entries of the second eigenvector $\mathbf{v_2}$ actually exceed that expectation by at least $c_{p,q}\cdot \sqrt{n} \cdot \eta$. Suppose not: then $1\%$ of the $\varepsilon-$largest entries of $\mathbf{v_2}$ actually have at most $c_{p,q}\cdot \sqrt{n} \cdot \eta$ neighbors; the value of their corresponding $\mathbf{v_2}$ entries must therefore originate in the error term, implying that
$$ \| \mathbf{error} \|^2 \gtrsim_{\eta, p, q} \frac{\varepsilon n}{100} \left( \frac{1}{n} \right)^2 \gtrsim \frac{1}{n}$$
which is a contradiction.  
\end{proof}


\begin{thebibliography}{10}

\bibitem{abbe} E. Abbe, 
Community detection and stochastic block models: recent developments.
J. Mach. Learn. Res. 18 (2017), Paper No. 177, 86 pp.


\bibitem{abbe2} E. Abbe, A. Bandeira, G. Hall, 
Exact recovery in the stochastic block model. 
IEEE Trans. Inform. Theory 62 (2016), no. 1, 471--487.



\bibitem{abbe1} E. Abbe, C. Sandon, 
Proof of the achievability conjectures for the general stochastic block model. 
Comm. Pure Appl. Math. 71 (2018), no. 7, 1334--1406.


\bibitem{band} A. Bandeira,
Random Laplacian matrices and convex relaxations.
Found. Comput. Math. 18 (2018), no. 2, 345--379.

\bibitem{blum} A. Blum, J. Hopcroft and R. Kannan, Foundations of Data Science, Cambridge University Press, 2020.

\bibitem{b25} K. Burdzy,  The hot spots problem in planar domains with one hole, Duke Math. J. 129 (2005),  p. 481--502.


\bibitem{b3}  R. Ba\~nuelos and K. Burdzy, On the "hot spots" conjecture of J. Rauch, J. Func. Anal. 164 (1999), p. 1--33



\bibitem{b4} K. Burdzy and W. Werner, A counterexample to the "hot spots" conjecture, Ann. Math. 149 (1999), p. 309--317

\bibitem{cheeger} J. Cheeger, A lower bound for the smallest eigenvalue of the Laplacian. Problems in analysis (Papers dedicated to Salomon Bochner, 1969), p. 195--199.

\bibitem{gal} X. Cheng, G. Mishne and S. Steinerberger, ]=The Geometry of Nodal Sets and Outlier Detection, arxiv, Journal of Number Theory, 185 , p. 48--64 (2018).

\bibitem{xiu} X. Cheng, M. Rachh and S. Steinerberger, On the Diffusion Geometry of Graph Laplacians and Applications, arxiv, Appl. Comp. Harm. Anal., 46, p. 674--688 (2019).

\bibitem{chung} F. R. K. Chung,  Spectral graph theory. CBMS Regional Conference Series in Mathematics, 92. Published for the Conference Board of the Mathematical Sciences, Washington, DC; by the American Mathematical Society, Providence, RI, 1997.

\bibitem{moo}  Moo Chung, Seongho Seo, Nagesh Adluru, and Houri Vorperian. Hot spots conjecture and its application to modeling tubular structures. In Kenji Suzuki, Fei
Wang, Dinggang Shen, and Pingkun Yan, editors, Machine Learning in Medical
Imaging, volume 7009 of Lecture Notes in Computer Science, p. 225--232.



\bibitem{fiedler1} 
M. Fiedler. Algebraic connectivity of graphs. Czechoslovak Math. J., 23(98), p. 298--305, 1973


\bibitem{fiedler2}
M. Fiedler. A property of eigenvectors of nonnegative symmetric matrices and its application to
graph theory. Czechoslovak Math. J., 25(100)(4), p.619--633, 1975.

\bibitem{fiedler3} 
M. Fiedler. Laplacian of graphs and algebraic connectivity. In Combinatorics and graph theory
(Warsaw, 1987), volume 25 of Banach Center Publ., p. 57--70. PWN, Warsaw, 1989.



\bibitem{gern} H. Gernandt and J. Pade, Schur reduction of trees and extremal entries of the
Fiedler vector, Linear Algebra and its Applications
Volume 570, p. 93--122 (2019)



\bibitem{ham} D. Hammond, P. Vandergheynst, R. Gribonval, Wavelets on graphs via spectral graph theory, Applied and Computational Harmonic Analysis 30, p. 129--150, (2011).

\bibitem{judge} C. Judge and S. Mondal, Euclidean Triangles Have No Hot Spots, Ann. of Math, to appear.

\bibitem{kato} T. Kato, Perturbation theory for linear operators.
Die Grundlehren der mathematischen Wissenschaften, Band 132 Springer-Verlag New York, Inc., New York 1966


\bibitem{kawohl} B. Kawohl, Rearrangements and Convexity of Level Sets in PDE, Lecture Notes in
Mathematics 1150, Springer, Berlin, 1985.

\bibitem{roy} R. Lederman and S. Steinerberger, Extreme Values of the Fiedler Vector on Trees, arXiv:1912.08327

\bibitem{levin}  D. Levin and Y. Peres, Markov chains and mixing times. With contributions by Elizabeth L. Wilmer. With a chapter on "Coupling from the past'' by James G. Propp and David B. Wilson. American Mathematical Society, Providence, RI, 2017.

\bibitem{lux} U. von Luxburg, A tutorial on spectral clustering, Statistics and Computing 17 (2007), p. 395--416.

\bibitem{mc} F. McSherry, Spectral partitioning of random graphs. (English summary) 42nd IEEE Symposium on Foundations of Computer Science (Las Vegas, NV, 2001), p. 529--537

\bibitem{ng} A. Ng, M. Jordan, Y. Weiss, On spectral clustering: analysis and an algorithm,  NIPS'01: Proceedings of the 14th International Conference on Neural Information Processing Systems: Natural and SyntheticJanuary 2001, p. 849--856 

\bibitem{manas} M. Rachh and S. Steinerberger, On the location of maxima of solutions of Schroedinger's equation, Comm. Pure. Appl. Math., 71, p.1109--1122 (2018).

\bibitem{rohe} K. Rohe, S. Chatterjee and Bin Yu, 
Spectral clustering and the high-dimensional stochastic blockmodel. 
Ann. Statist. 39 (2011), no. 4, p. 1878--1915.

\bibitem{spiel} D. Spielman, S.-H. Teng, Spectral partitioning works: planar graphs and finite element meshes. 37th Annual Symposium on Foundations of Computer Science (Burlington, VT, 1996), p. 96--105, IEEE Comput. Soc. Press, Los Alamitos, CA, 1996.

\bibitem{steini} S. Steinerberger, Hot Spots in Convex Domains are in the Tips (up to an Inradius), arXiv:1907.13044

\bibitem{vershynin} R. Vershynin, High-dimensional probability. An introduction with applications in data science. Cambridge Series in Statistical and Probabilistic Mathematics, 47. Cambridge University Press, Cambridge, 2018.
\end{thebibliography}
\end{document}